\documentclass[12pt]{article}
\usepackage[numbers]{natbib}
\usepackage{graphicx}
\usepackage{color, multirow,amsfonts,amsmath}
\usepackage{url}
\usepackage{rotating}
\title{Mixed Poisson process with Min-U-Exp mixing variable - Work in progress
\bigskip
\\Pavlina K. Jordanova
\\{\small{\it Faculty of Mathematics and Informatics, Konstantin Preslavsky University of Shumen, \\115 "Universitetska" str., 9712 Shumen, Bulgaria. \\Corresponding author:  pavlina\_kj@abv.bg.}}
\\Evelina Veleva
\\{\small{\it Department of Applied mathematics and Statistics, Angel Kanchev  University of Ruse, Bulgaria.}}
\\Milan Stehlik
\\{\small{\it Institute of Applied Statistics, Johannes Kepler University in Linz, Austria.}}
}

\begin{document}
\date{}

\maketitle

\begin{abstract}
This work continues the research done in Jordanova and Veleva (2023) where the history of the problem could be found. In order to obtain the structure distribution of the newly-defined Mixed Poisson process, here the operation "max" is replaced with "min". We start with the definition of Min-U-Exp distribution. Then, we compute its numerical characteristics and investigate some of its properties. The joint distribution of the inter-arrival times (which are dependent) is the Multivariate  Exp-Min-U-Exp distribution of $II^{-nd}$ kind. Its univariate and multivariate versions are described, and the formulae for their numerical characteristics are obtained. The distribution of the moments of arrival of different events is called Erlang-Min-U-Exp.  Different properties of these distributions are obtained, and their numerical characteristics are computed. Multivariate ordered Mixed Poisson-Min-U-Exp distribution describes the joint distribution of the time-intersection of a Mixed Poisson process with Min-U-Exp mixing variable. The corresponding distribution of the additive increments (which are also dependent) is the Mixed Poisson-Min-U-Exp one.  The considered relations between these distributions simplify their understanding.
\end{abstract}

\section{DESCRIPTION OF THE MODEL AND PRELIMINARIES}
Let $\mathfrak{u}$ be a Uniformly distributed random variable (r.v.) on the interval $(0, a)$, briefly $\mathfrak{u} \in U(0, a)$. Here and further on, we denote by $\mathfrak{e}$ an Exponentially distributed r.v. with mean $\frac{1}{\lambda}$, $\lambda > 0$, i.e. $\mathfrak{e} \in Exp(\lambda)$. Analogously to Jordanova and Veleva (2023) \cite{JV2023}, by replacing the maxima with minima, we consider the distribution of $\xi := \min(\mathfrak{u}, \mathfrak{e})$, and we call it Min-U-Exp distribution.

{\bf Definition 1.} \label{Def:1} We say that the r.v. $\xi$ is {\bf{Min-U-Exp distributed}} with parameters $a > 0$ and $\lambda > 0$, if it has a cumulative distribution function (c.d.f.)
\begin{equation}\label{MinUExpCDF}
F_\xi(x) = \left\{ \begin{array}{ccc}
                                   0 & , & x \leq 0\\
                                   1-e^{-\lambda x} + \frac{x}{a}e^{-\lambda x}& , & x \in (0, a]\\
                                   1& , & x > a\\
                      \end{array}\right..
\end{equation}

Briefly, we will denote this in this way  $\xi \in Min-U-Exp(a; \lambda)$.

In the next proposition and theorem we investigate the main properties of Min-U-Exp distribution. The proves are analogous to the corresponding one in Jordanova and Veleva (2023).

{\bf Proposition 1.}
\begin{description}
  \item[a)] $\xi \in Min-U-Exp(a; \lambda)$ if and only if the probability density function (p.d.f.) of $\xi$ is
\begin{equation}\label{MinUExpDensity}
P_\xi(x) = \left\{ \begin{array}{ccc}
                                   0 & , & x \not\in (0,a)\\
                                   \frac{e^{-\lambda x}}{a}(\lambda a + 1 - x \lambda)& , & x \in (0, a)
           \end{array}\right..
\end{equation}
  \item[b)] (Scaling property) If $\xi \in Min-U-Exp(a; \lambda)$ and $k > 0$ is a constant, then
  $$k\xi \in Min-U-Exp\left(ka; \frac{\lambda}{k}\right).$$
  \item[c)]  If $\xi \in Min-U-Exp(a; \lambda)$, the hazard rate function of this distribution is
 $$h_{\xi}(x) = \left\{\begin{array}{ccc}
                   0 & , & x \not\in (0,a)\\
                   \lambda + \frac{1}{a-x} & , & x \in (0, a]
                    \end{array}
  \right..$$
\end{description}

\begin{minipage}{8cm}
\begin{center}
    \includegraphics[scale=.65]{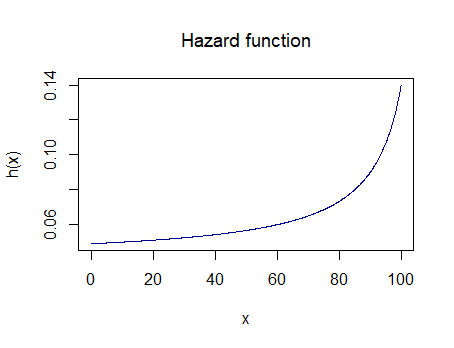}

 \small{Figure 1. Hazard function of $\xi \in Min-U-$ $Exp(110; 0.04)$. \label{fig:Fig1HazardFunctionMinUExp}}
\end{center}
\end{minipage} \hfill
\begin{minipage}{7.5cm}
 {\it Note:} The hazard function described in c), and plotted in Figure 1, for $\lambda = 0.04$ and $a = 110$, is increasing and $$\lim_{x \uparrow a} h_{\xi}(x) = \infty.$$ It is well-known, that in the terms of survival theory this means that if $\xi$ describes the length of someone's life, then, the risk for his/her death is higher as his/her age increases, and when the age approaches $a$, this risk becomes infinite. The fact that $$\lim_{x \downarrow 0} h_{\xi}(x) = \lambda + \frac{1}{a}$$ means that in the beginning of his/her life also there is some possibility for death. 
 \end{minipage}
\smallskip

These conclusions allow us to say that, for different parameters, this distribution is very appropriate for modelling of many lengths of lives.
\bigskip

{\bf Theorem 1.} Let $a > 0$, $\lambda > 0$, $\mathfrak{u} \in U(0, a)$, $\mathfrak{e} \in Exp(\lambda)$, and $\mathfrak{u}$ and $\mathfrak{e}$ be independent r.vs.
 Denote by $\xi := \min(\mathfrak{u}, \mathfrak{e})$. Then,
\begin{itemize}
\item [a)] $\xi \in Min-U-Exp(a; \lambda)$;
\item [b)] The mean, and the moments of $\xi$ are correspondingly  $\mathbb{E}\xi = \frac{1}{a \lambda^2}(a\lambda - 1 + e^{-\lambda a}),$ and
$$\mathbb{E}(\xi^k) = \frac{k}{\lambda^k}\left(\gamma(k, a\lambda) - \frac{\gamma(k+1, a\lambda)}{a\lambda}\right), \quad k \in \mathbb{N}.$$
\item [c)] The variance of $\xi$ is $\mathbb{D}\xi = \frac{1}{\lambda^2}\left(2 + 2e^{-\lambda a}-\frac{1}{a^2\lambda^2}(a\lambda+1-e^{-\lambda a})^2\right) .$
\item [d)] The Laplace-Stieltjes transform of $\xi$ is
  $\mathbb{E}(e^{-\xi t}) = \frac{\lambda}{\lambda + t} + \frac{t}{a(\lambda + t)^2}(1 - e^{-(\lambda + t) a})$, $t \geq 0$.
\end{itemize}

{\bf Proof:} a) Consider $x \in \mathbb{R}$, the definition of $\xi$ and the independence between $\mathfrak{u}$ and $\mathfrak{e}$ entail,
$F_\xi(x) =  \mathbb{P}(\min(\mathfrak{u}, \mathfrak{e}) \leq x) = \mathbb{P}(\mathfrak{u} \leq x \cup \mathfrak{e} \leq x)  = 1 - \mathbb{P}(\mathfrak{u} > x \cap \mathfrak{e} > x) = 1 - \mathbb{P}(\mathfrak{u} > x)\mathbb{P}(\mathfrak{e} > x) = 1 - (1-\mathbb{P}(\mathfrak{u} \leq x))(1- \mathbb{P}(\mathfrak{e} \leq x)).$
The definitions of $Exp(\lambda)$ and $U(0, a)$ distributions via their c.d.fs. entail (\ref{MinUExpCDF}) and complete the proof of a) .

b) {\footnote{Although this proof could be analogous to the corresponding one in Jordanova and Veleva (2023) \cite{JV2023} here we present a different approach.}
By taking expectation in the both sides of the equality $\min(\mathfrak{u}, \mathfrak{e}) + \max(\mathfrak{u},\mathfrak{e}) = \mathfrak{u} + \mathfrak{e}$, and by using the additive property of the expectations we obtain, $\mathbb{E}(\min(\mathfrak{u}, \mathfrak{e})) + \mathbb{E}(\max(\mathfrak{u},\mathfrak{e})) = \mathbb{E}\mathfrak{u} + \mathbb{E}\mathfrak{e}$. Now, the Theorem 1, b) in Jordanova and Veleva (2023) \cite{JV2023}, and the well-known formulae for the expectations of the exponential and uniform distributions lead us to the equality $\mathbb{E}(\min(\mathfrak{u}, \mathfrak{e})) + \frac{a}{2} + \frac{1}{a \lambda^2}(1 - e^{-\lambda a}) = \frac{a}{2} + \frac{1}{\lambda}$.

For all $k \in \mathbb{N}$, the equality $\min(\mathfrak{u}, \mathfrak{e})^k + \max(\mathfrak{u},\mathfrak{e})^k = \mathfrak{u}^k + \mathfrak{e}^k$, entails $\mathbb{E}(\min(\mathfrak{u}, \mathfrak{e})^k) + \mathbb{E}(\max(\mathfrak{u},\mathfrak{e})^k) = \mathbb{E}(\mathfrak{u}^k) + \mathbb{E}(\mathfrak{e}^k)$. Analogously to the proof of the expectations $\mathbb{E}(\min(\mathfrak{u}, \mathfrak{e})^k) + \frac{a^k}{k+1} + \frac{k}{a\lambda^{k+1}}\gamma(k+1, a\lambda) + \frac{k}{\lambda^k}\Gamma(k, \lambda a)= \frac{a^k}{k+1} + \frac{k!}{\lambda^k}$. The rest follows by the well-known relation $k! = k\Gamma(k) = k\gamma(k, a\lambda) + k\Gamma(k, a\lambda)$.

c) After some algebra, the relation $\mathbb{D}\xi = \mathbb{E}(\xi^2) - (\mathbb{E}\xi)^2$ and b) entail c).

d$)^1$ For all $t \geq 0$, the equality $e^{-t \min(\mathfrak{u}, \mathfrak{e})} + e^{-t \max(\mathfrak{u},\mathfrak{e})} = e^{-t\mathfrak{u}} + e^{-t\mathfrak{e}}$, and the additive property of the mean entail
$$\mathbb{E}e^{-t \min(\mathfrak{u}, \mathfrak{e})} + \mathbb{E}e^{-t \max(\mathfrak{u},\mathfrak{e})} = \mathbb{E}e^{-t\mathfrak{u}} + \mathbb{E}e^{-t\mathfrak{e}}.$$
Now, the Theorem 1, d) in Jordanova and Veleva (2023) \cite{JV2023}, and the well-known formulae for the Laplace-Stieltjes transform of the exponential and uniform distributions lead us to
\begin{eqnarray*}
\mathbb{E}(e^{-t\xi}) & = &\mathbb{E}(e^{-t\min(\mathfrak{u}, \mathfrak{e})}) =  \frac{e^{-at}-1}{ta} + \frac{\lambda}{\lambda + t} - \frac{1}{at}(1 - e^{-\lambda a}) +\frac{t}{a(\lambda + t)^2}(1 - e^{-(\lambda + t) a})\\
& = & \frac{\lambda}{\lambda + t} +\frac{t}{a(\lambda + t)^2}(1 - e^{-(\lambda + t) a}).
\end{eqnarray*}
  \hfill $\Box$

\bigskip

The method of moments can be used to estimate the parameters of $Min-U-Exp(a; \lambda)$ distribution. Suppose we have a sample of size $n$ of observations over this distribution. The first two empirical initial moments $m_1$ and $m_2$ are unbiased and consistent estimates of $E\xi$ and $E\xi^2$, the first and second initial moments of the distribution. According to Theorem 1, $E\xi$ and $E\xi^2$ are equal to respectively

\begin{equation}\label{3}
\left|\begin{matrix}E\xi=\frac{1}{a\lambda^2}\left(a\lambda-1+e^{-a\lambda}\right)=\frac{1}{xy}\left(x-1+e^{-x}\right)\\E{(\xi}^2)=\frac{2}{a\lambda^3}\left(a\lambda-2+a\lambda e^{-a\lambda}+2e^{-a\lambda}\right)=\frac{2}{xy^2}\left(x-2+{\rm xe}^{-x}+{2e}^{-x}\right)\\\end{matrix}\right.
\end{equation}
where  $x = a\lambda$ and $y = \lambda$. When $E\xi$ and $E\xi^2$ are replaced by their empirical estimates $m_1$ and $m_2$, a nonlinear system of two equations with two unknowns $a$ and $\lambda$ (or $x$ and $y$) is obtained, which can be easily solved using, for example, \textit{Matlab} or \textit{Octave} software and the \textit{fsolve} command. To investigate whether system (\ref{3}) will have one or more solutions, let us write it in equivalent form:
\begin{equation}\label{4}
\left|\begin{matrix}2x\left(x-2+{\rm xe}^{-x}+{2e}^{-x}\right):\left(x-1+e^{-x}\right)^2=\frac{E{(\xi}^2)}{{(E\xi)}^2}\\y=\left(x-1+e^{-x}\right):(xE\xi)\\\end{matrix}\right.
\end{equation}

The left side of the first equation of system (\ref{4}) is a non-linear function, say $G(x)$, of $x = a\lambda$. Its graph is given in Figure 2. By Taylor expansion series about zero, it can be established that $\lim_{x \downarrow 0} G(x) = 4/3 = 1.3333$. On the other hand, since $G(x)$ is quotient to two polynomials of equal degrees, $\lim_{x\rightarrow\infty} G(x)=2$. Therefore, the right-hand side of the first equation, equal to $\frac{E{(\xi}^2)}{{(E\xi)}^2}$, will always be a number in the interval $(4/3, 2)$. To each estimate $\hat{r}$ of the ratio $\frac{E{(\xi}^2)}{{(E\xi)}^2}$, $\hat{r} \in (4/3, 2)$, there will correspond a unique value $x^*$, a solution of the equation $G(x)=\hat{r}$. The solution $x^*$ can be found numerically (for example with the \textit{fzero} command of \textit{Matlab} and \textit{Octave}) or even graphically from Figure 2. Then, from the second equation of system (\ref{4}) we find the corresponding value of $y = \lambda$, replacing in the right hand side $x$ with $x^*$ and $E\xi$ with $m_1$. Finally, we determine the value for the parameter $a$: $a=x^*/y=(a\lambda)/\lambda$. When the estimator $\hat{r} \in (4/3, 2)$, the system (\ref{4}) will have an unique solution. For $\hat{r} \le 4/3$ it can be assumed that $\lambda  =0$ and $a = \xi^{(n)}$, the maximum observation in the sample, since $\xi^{(n)}$ is the maximum likelihood estimate for the right end of the definition interval in the uniform distribution.

\begin{center}
\includegraphics[scale=.56]{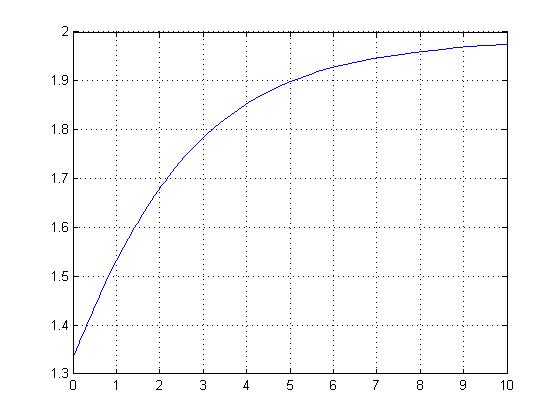}

\small{Figure 3. Graph of the left-hand side of the first equation of system (\ref{4}) as a function $G(x)$ of $x = a\lambda$ on the interval (0, 10] \label{fig:fig2Eva}} 
\end{center}

Alternatively, the method of least squares can be used to estimate the parameters of $Min-U-Exp(a; \lambda)$ distribution. For each of the $n$ observations in the sample, the value of the empirical distribution function is calculated. For estimates of parameters $a$ and $\lambda$, the values minimizing the sum of the squares of the differences between the values of the empirical and theoretical distribution functions (corresponding to the observations in the sample) are taken. This approach can also be easily implemented using software such as \textit{Matlab} and command \textit{fit}.

\section{EXP-MIN-U-EXP AND ERLANG-MIN-U-EXP DISTRIBUTIONS}

{\bf Definition 2.} \label{Def:2} We say that the r.v. $\tau$ is {\bf{Exp-Min-U-Exp distributed}} with parameters $a > 0$ and $\lambda > 0$, if it has a c.d.f.
\begin{equation}\label{ExpMinUExpCDF}
F_\tau(t) = \left\{ \begin{array}{ccc}
                                   0 & , & t \leq 0\\
                                   \frac{t}{\lambda+ t} - \frac{t}{a(\lambda + t)^2}(1-e^{-a(\lambda + t)})& , & t > 0
                      \end{array}\right..
\end{equation}

Briefly we will denote this in this way  $\tau \in Exp-Min-U-Exp(a; \lambda)$.

The well-known relations between c.d.fs., p.d.fs. and the corresponding  probability distribution entail the following result.

{\bf Proposition 2.} For $a > 0$ and $\lambda > 0$, $\tau \in Exp-Min-U-Exp(a; \lambda)$ if and only if the p.d.f.
\begin{equation}\label{ExpMinUExpDensity}
P_\tau(t) = \left\{ \begin{array}{ccc}
                                   0 & , & t \leq 0\\
                                   \frac{\lambda}{(\lambda + t)^2}+\frac{t-\lambda}{a(\lambda + t)^3}(1-e^{-a(\lambda + t)})-\frac{t}{(\lambda + t)^2}e^{-a(\lambda + t)}& , & t > 0
                      \end{array}\right..
\end{equation}

Analogously to the corresponding result in \cite{JV2023} we obtain that $\int_0^\infty P_\tau(t) dt = 1$. The last means that this distribution is proper.

{\bf Definition 3.} \label{Def:3} We say that the random vector (rv.) $(\tau, \xi)$ has {\bf{bivariate Exp-Min-U-Exp distribution of $I^{-st}$ kind}} with parameters $a > 0$, and $\lambda > 0$, if it has a joint p.d.f.
\begin{equation}\label{BivariateExpMaxUExpDensity}
P_{\tau, \xi}(t,x) = \left\{ \begin{array}{ccc}
                                   \frac{1}{a}xe^{-(\lambda + t)x}(1+\lambda a- x \lambda) & , & t > 0 \cap x \in (0;a)\\
                                   0 & , & otherwise\\
                      \end{array}\right..
\end{equation}
Briefly we will denote this in this way  $(\tau, \xi) \in Exp-Max-U-Exp-I^{-st}(a, \lambda)$.

The proves of the results in the next theorem are analogous to the corresponding one in \cite{JV2023}. Here we are going to present only some different approaches for some of them.

{\bf Theorem 2.} For $a > 0$ and $\lambda > 0$, if $\xi \in Min-U-Exp(a; \lambda)$ and for $x > 0$, $(\tau|\xi=x) \in Exp(x)$, then:
\begin{description}
  \item[a)] $\tau \in Exp-Min-U-Exp(a, \lambda)$;
  \item[b)] $\tau \stackrel{d}{=} \frac{\eta}{\xi}$, where $\eta \in Exp(1)$, and $\xi$ and $\eta$ are independent.
  \item[c)] For $p \in (-1, 1)$,
  $$\mathbb{E}(\tau^p) = \frac{1}{a}\Gamma(p+1)\lambda^{p-1}\left((\lambda a+1)\gamma(1-p, a\lambda) -\gamma(2-p, a\lambda)\right),$$
   and $\mathbb{E}(\tau^p) = \infty$, otherwise.
  \item[d)] The joint distribution of $\tau$ and $\xi$ is $(\tau, \xi) \in Exp-Min-U-Exp-I^{-st}(a, \lambda)$ and $(\tau, \xi) \stackrel{d}{=} \left(\frac{\eta}{\xi}, \xi\right)$, where $\eta \in Exp(1)$, and $\xi$ and $\eta$ are independent.
  \item[e)] For all $t > 0$,
    $$P_{\xi}(x|\tau=t) = 0, \quad x \not\in (0;a),$$
    $$P_{\xi}(x|\tau=t) = \frac{x(\lambda + t)^3e^{-(\lambda + t)x}(1+\lambda a - \lambda x)}{a\lambda(\lambda + t) + (t-\lambda)(1-e^{-a(\lambda +t)})- a t (\lambda+t)e^{-a(\lambda + t)}}, \quad  x \in (0, a).$$
   \item[f)] The mean square regression $\mathbb{E}(\tau|\xi=x) = \frac{1}{x}$, $x > 0$.
   \item[g)] The mean square regression function is
   $$\mathbb{E}(\xi|\tau=t) = \frac{2(t-2\lambda) +2a\lambda(t+\lambda)-e^{-a(\lambda +t)}(a^2t(\lambda+t)^2+2a(t^2-\lambda^2)-2(t-2\lambda))}{a\lambda(\lambda+t)^2+(t^2-\lambda^2)(1-e^{-a(\lambda + t)})-at(\lambda+t)^2e^{-a(\lambda + t)}}, \quad t > 0.$$
\end{description}

{\bf Proof:} a) For $t > 0$, the integral form of the Total probability formula and Theorem 1, d) entail,
\begin{eqnarray*}
  \mathbb{P}(\tau > t) &=& \int_0^\infty \mathbb{P}(\tau > t|\xi=x)P_\xi(x)dx = \int_0^\infty e^{-\lambda x}P_\xi(x)dx = \mathbb{E}(e^{-\xi t})\\
   &=& \frac{\lambda}{\lambda + t} + \frac{t}{a(\lambda + t)^2}(1 - e^{-(\lambda + t) a}).
\end{eqnarray*}
Thus, the relation $F_\tau(t) = 1- \mathbb{P}(\tau > t)$ leads us to (\ref{ExpMinUExpCDF}), and $\tau \in Exp-Min-U-Exp(a, \lambda)$.  \hfill $\Box$

 {\bf Definition 4.} \label{Def:4} We say that a rv. $(\tau_1, \tau_2, ..., \tau_k)$ has {\bf Multivatiate Exp-Min-U-Exp distribution of $II^{-nd}$ kind with parameters $a > 0$, and $\lambda > 0$}, if it has a joint p.d.f.

$$P_{\tau_1, \tau_2, \ldots, \tau_k}(t_1, t_2, \ldots, t_k) = \gamma(k+1,a(\lambda + t_1 + \ldots + t_k))\frac{a\lambda(\lambda + t_1 + \ldots + t_k) + t_1 + \ldots + t_k - \lambda k}{a(\lambda + t_1 + \ldots + t_k)^{k+2}} $$
$$+ \frac{\lambda a^k}{\lambda +  t_1 + \ldots + t_k}e^{-a(\lambda + t_1 + \ldots + t_k)},
\quad t_1 > 0, t_2 > 0, \ldots, t_k > 0,$$
and $P_{\tau_1, \tau_2, \ldots, \tau_k}(t_1, t_2, \ldots, t_k) = 0$,  otherwise.

Briefly we will denote this in this way  $(\tau_1, \tau_2, \ldots, \tau_k) \in Exp-Min-U-Exp-II(a, \lambda)$.

{\bf Definition 5.} \label{Def:5} We say that the r.v. $T_n$ is {\bf{Erlang-Min-U-Exp distributed with parameters $n \in \mathbb{N}$, $a > 0$, and $\lambda > 0$}}, if it has a p.d.f.
$$P_{T_n}(t) = \frac{t^{n-1}\gamma(n+1,a(\lambda + t))}{a(n-1)!(\lambda+t)^{n-1}}\left(\lambda a + \frac{t - \lambda n}{\lambda+t}\right) + \frac{\lambda a^n t^{n-1}}{(n-1)!(\lambda + t)}e^{-a(\lambda + t)}, $$
when $t > 0$, and $P_{T_n}(t) = 0$, otherwise. Briefly, we will denote this in this way  $T_n \in Erlang-Min-U-Exp(n; a, \lambda)$.

We will skip the proves of the results in the next theorem as far as they are analogous to the corresponding one in \cite{JV2023}.

{\bf Theorem 3.} For $a > 0$, and $\lambda > 0$, if $\xi \in Min-U-Exp(a; \lambda)$ and for $x > 0$, $(\tau_1, \tau_2, ..., \tau_k|\xi=x)$ are independent identically $Exp(x)$ distributed r.vs.,  $T_n := \tau_1 + \ldots + \tau_n$, $n \in \mathbb{N}$, then,

\begin{description}
  \item[a)] $(\tau_1, \tau_2, \ldots, \tau_k) \in Exp-Min-U-Exp-II(a, \lambda)$.
  \item[b)] For all $i = 1, 2, ..., k$, $\tau_i \in  Exp-Min-U-Exp-(a, \lambda)$.
  \item[c)]  $(\tau_1, \tau_2, \ldots, \tau_k)  \stackrel{d}{=} \left(\frac{\eta_1}{\xi}, \frac{\eta_2}{\xi}, \ldots, \frac{\eta_k}{\xi}\right)$, where $\eta_1, \eta_2, \ldots, \eta_k$ are independent identically distributed (i.i.d.) $Exp(1)$, and independent on $\xi$.
 \item[d)] $T_n \in Erlang-Min-U-Exp(n; a, \lambda)$. $T_n \stackrel{d}{=} \frac{\eta_1 + \eta_2 + \ldots + \eta_n}{\xi}$, where $\eta_1, \eta_2, \ldots, \eta_n$ are i.i.d. $Exp(1)$, and independent on $\xi$. $T_n \stackrel{d}{=} \frac{\theta_n}{\xi}$, where $\theta_n \in Gamma(n, 1)$ is independent on $\xi$.
   \item[e)] For $p \in (-n, 1)$,
 $$ \mathbb{E}(T_n^p) = \frac{\Gamma(p+n)}{(n-1)!}\lambda^p\left\{\left(1 + \frac{p}{a\lambda}\right)\gamma(1-p,\lambda a)+\frac{e^{-\lambda a}}{(a\lambda)^p}\right\},$$
  and $\mathbb{E}(T_n^p) = \infty$, for $p \not\in (-n, 1)$.
  \end{description}

\bigskip
\section{THE MIXED POISSON-MIN-U-EXP PROCESS}

{\bf Definition 6.} \label{Def:6} A r.v. $\theta$ has a {\bf{Mixed Poisson-Min-U-Exp distributed with parameters $a > 0$, and $\lambda > 0$}} if its probability mass function (p.m.f.) is
\begin{equation}\label{MixedPoissonUEXP}
\mathbb{P}(\theta = n) = \frac{1}{n!}\left\{\frac{\gamma(n+1,a(\lambda+1))}{(\lambda+1)^{n+2}}\left(\lambda(\lambda+1)+\frac{1-n\lambda}{a}\right) + \frac{\lambda a^n}{\lambda+1}e^{-a(\lambda+1)}\right\}, \quad n = 0, 1, \ldots.
\end{equation}
Briefly, $\theta \in MPMin-U-Exp(a, \lambda)$.

\begin{minipage}{8cm}
\begin{center}
    \includegraphics[scale=.8]{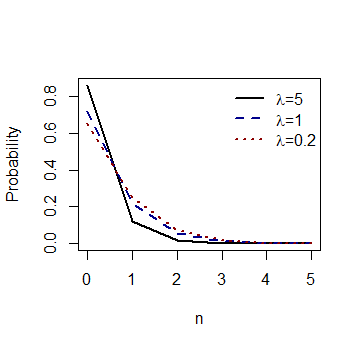}

 \small{Figure 4. P.m.f. of $\theta \in MPMin-U-Exp(1, \lambda)$. \label{fig:Fig2MPPFixeda1}}
\end{center}
\end{minipage} \hfill
\begin{minipage}{8cm}
\begin{center}
    \includegraphics[scale=.8]{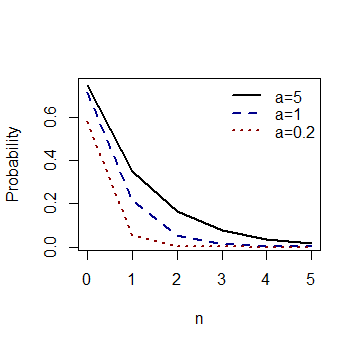}

 \small{Figure 5. P.m.f. of $\theta \in MPMin-U-Exp(a, 1)$. \label{fig:Fig3MPPFixedlambda1}}
\end{center}
\end{minipage}

\bigskip

{\bf Definition 7.} \label{Def:7} Let $\mu(t): [0, \infty) \to [0, \infty)$ be a nonnegative, strictly increasing and continuous function, $\mu(0) = 0$, $\xi \in Min-U-Exp(a; \lambda)$ and $N_1$  be a Homogeneous Poisson process (HPP) with intensity $1$, independent on $\xi$. We call the random process
\begin{equation}\label{TheProcessN}
N := \{N(t), t\geq 0\}  = \{N_1(\xi \mu(t)), t \geq 0\}
\end{equation}
a {\bf{Mixed Poisson process with Min-U-Exp mixing variable}} or {\bf{MPMin-U-Exp process}}. Briefly $N \in MPMin-U-Exp(a, \lambda; \mu(t))$.

\medskip

{\bf Definition 8.} \label{Def:8} Let $n \in \mathbb{N}$. We say that a random vector $(N_1, N_2, \ldots, N_n)$ is {\bf Ordered Poisson-Min-U-Exp distributed with parameters $a > 0$, $\lambda > 0$, and $0 < \mu_1 < \mu_2 < ... < \mu_n$} if, for all integers $0 \leq k_1 \leq k_2 \leq \ldots \leq k_n$,
$$\mathbb{P}(N_1 = k_1, N_2 = k_2, \ldots, N_n = k_n) = \frac{\mu_1^{k_1}(\mu_2 - \mu_1)^{k_2 - k_1}\ldots(\mu_n - \mu_{n-1})^{k_n - k_{n - 1}}}{k_1!(k_2 - k_1)!\ldots(k_n - k_{n-1})!}$$
$$\times \left\{\frac{\gamma(k_n+1,a(\mu_n+\lambda))}{a(\mu_n+\lambda)^{k_n+2}}(\lambda a(\lambda+\mu_n) + \mu_n-\lambda k_n)+\frac{a\lambda}{\lambda+\mu_n}e^{-a(\lambda+\mu_n)}\right\},$$
and $\mathbb{P}(N_1 = k_1, N_2 = k_2, \ldots, N_n = k_n) = 0$, otherwise.

Briefly, $(N_1, N_2, \ldots, N_n) \in O_{PMinUE}(a, \lambda; \mu_1, \mu_2, ..., \mu_n)$.

\medskip

{\bf Definition 9.} \label{Def:9} Let $n \in \mathbb{N}$. We say that a random vector $(N_1, N_2, \ldots, N_n)$ is {\bf Mixed Poisson-Min-U-Exp  distributed with parameters $a > 0$, $\lambda > 0$, and $0 < \mu_1 < \mu_2 < ... < \mu_n$} if, for all $m_1, m_2, \ldots, m_n \in \{0, 1, \ldots\}$,
$$\mathbb{P}(N_1 = m_1, N_2 = m_2, \ldots, N_n = m_n) = \frac{\mu_1^{m_1}(\mu_2 - \mu_1)^{m_2}\ldots(\mu_n - \mu_{n-1})^{m_n}}{ m_1!m_2!\ldots m_n!}$$
$$\times \left\{\frac{\gamma(m_1 + \ldots + m_n + 1,a(\lambda+\mu_n))}{a(\lambda + \mu_n)^{m_1 + \ldots + m_n + 2}}(\lambda a(\lambda+\mu_n) + \mu_n-\lambda(m_1 + \ldots + m_n)) + \frac{a  \lambda}{\lambda + \mu_n}e^{-a(\lambda + m_n)}\right\},$$
and $\mathbb{P}(N_1 = m_1, N_2 = m_2, \ldots, N_n = m_n) = 0$, otherwise.

Briefly, $(N_1, N_2, \ldots, N_n) \in M_{PMinUE}(a, \lambda; \mu_1, \mu_2, ..., \mu_n)$.

\medskip
 In the next two propositions we present two relations between the distributions introduced in  Definition 8 and Definition 9. Their proves, together with the proof of Theorem 4 will be skipped, because they are analogous to the corresponding one in Jordanova et al. \cite{SMladen}, and Jordanova and Stehlik \cite{LMJ}. The algorithms are based on the above results and the general formulae for any Mixed Poisson process which could be found, for example, in Grandel \cite{GrandelMixed}, or Karlis and Xekalaki \cite{KarlisandXekalaki}.

{\bf Proposition 3.} If $(N_1, N_2, \ldots, N_n) \in O_{PMinUE}(a, \lambda; \mu_1, \mu_2, ..., \mu_n)$, then
$$(N_1, N_2 - N_1, \ldots, N_n - N_{n-1}) \in M_{PMinUE}(a, \lambda; \mu_1, \mu_2, ..., \mu_n).$$

{\bf Proposition 4.}  If $(N_1, N_2, \ldots, N_n) \in M_{PMinUE}(a, \lambda; \mu_1, \mu_2, ..., \mu_n)$, then
$$(N_1, N_1 + N_2, \ldots, N_1 + N_2 + \ldots + N_n) \in O_{PMinUE}(a, \lambda; \mu_1, \mu_2, ..., \mu_n).$$

{\bf Theorem 4.} Consider $a > 0$, $\lambda > 0$, and  a nonnegative, strictly increasing and continuous deterministic function $\mu(t): [0, \infty) \to [0, \infty)$. Suppose that $\{N(t), t \geq 0\} \in MPMin-U-Exp(a, \lambda; \mu(t))$.
\begin{description}
\item[a)] For all $t > 0$, $N(t) \in MPMin-U-Exp(a\mu(t), \frac{\lambda}{\mu(t)})$.
\item[b)]  These processes are over-dispersed,
$$\mathbb{E}N(t) = \frac{\mu(t)}{\lambda}\left(1 - \frac{1}{a\lambda} + \frac{e^{-\lambda a}}{a\lambda}\right),$$
$$\mathbb{D}N(t) = \frac{\mu(t)}{\lambda}\left(1 - \frac{1}{a\lambda} + \frac{e^{-\lambda a}}{a\lambda}\right) + \frac{\mu^2(t)}{\lambda^2}\left\{2+2e^{-\lambda a}-\left(1+ \frac{1}{a\lambda} - \frac{e^{-\lambda a}}{a\lambda}\right)^2 \right\}.$$
\item[c)] For all $t \geq 0$, the probability generating function (p.g.f.) of the time intersections is
$$\mathbb{E}(z^{N(t)}) = \frac{\lambda}{\lambda + \mu(t)(1-z)} + \frac{\mu(t)(1-z)}{a(\lambda + \mu(t)(1-z))^2}\left(1 - e^{-a(\lambda + \mu(t)(1-z))}\right), \quad |z| < 1.$$
\item[d)] For $t > 0$, and $n = 0, 1, \ldots$, $P_{\xi}(x|N(t) = n) = 0$, when $x \leq 0$ or $x > a$, and when $x \in (0, a]$,
 $$P_{\xi}(x|N(t) = n) = \frac{x^ne^{-x(\lambda + \mu(t))}(\lambda a + 1 - \lambda x)}{\frac{\gamma(n+1,a(\lambda + \mu(t)))}{(\lambda + \mu(t))^{n+2}}\left(a \lambda(\lambda+\mu(t))+\mu(t) - n\lambda\right) + \frac{\lambda a^{n+1}}{\lambda+\mu(t)}e^{-a(\lambda+\mu(t))}}.$$
\item[e)] For $t > 0$, and $n = 0, 1, \ldots$, the mean square regression is
$$\mathbb{E}(\xi|N(t) = n) = \frac{\frac{\gamma(n+2,a(\lambda + \mu(t)))}{(\lambda + \mu(t))^{n+3}}\left(a \lambda(\lambda+\mu(t))+\mu(t) - (n+1)n\lambda\right) + \frac{\lambda a^{n+2}}{\lambda+\mu(t)}e^{-a(\lambda+\mu(t))}}{\frac{\gamma(n+1,a(\lambda + \mu(t)))}{(\lambda + \mu(t))^{n+2}}\left(a \lambda(\lambda+\mu(t))+\mu(t) - n\lambda\right) + \frac{\lambda a^{n+1}}{\lambda+\mu(t)}e^{-a(\lambda+\mu(t))}}.$$
\item[f)]  For all $k = 1, 2, \ldots$,
$$\mathbb{E}[N(t)(N(t)-1)(N(t) - k + 1)] = \frac{k(\mu(t))^k}{\lambda^k} \left(\gamma(k, a\lambda) - \frac{\gamma(k+1, a\lambda)}{a\lambda}\right).$$
\item[g)]  For all $n \in \mathbb{N}$, and $0 \leq t_1 \leq t_2 \leq \ldots \leq t_n,$
$$(N(t_1), N(t_2), \ldots, N(t_n)) \in O_{PMinUE}(a, \lambda; \mu(t_1), \mu(t_2), ..., \mu(t_n)).$$
\item[h)] For all $n \in \mathbb{N}$, and $0 \leq t_1 \leq t_2 \leq \ldots \leq t_n$,
$$(N(t_1), N(t_2) - N(t_1), \ldots, N(t_n) - N(t_{n-1})) \in M_{PMinUE}(a, \lambda; \mu(t_1), \mu(t_2), ..., \mu(t_n)).$$
\item[i)] Denote by $\tau_1, \tau_2, \ldots$ the inter-occurrence times of the counting process $N$. Then, $\tau_1, \tau_2, \ldots$ are dependent and $Exp-Min-U-Exp(a; \lambda)$ distributed.
\item[j)]  For $n \in \mathbb{N}$, if $T_n$ is the moment of occurrence of the $n$-th event of the counting process $N$, then $T_n \in Erlang-Min-U-Exp(n; a, \lambda)$.
\item[k)]  For any $0 < s < t$, the r.v. $N(s)$ given $N(t) = n$ is Binomially distributed. More precisely,
$$(N(s)|N(t) = n) \in Bi\left(n, \frac{\mu(s)}{\mu(t)}\right).$$
\end{description}

\bigskip
\section{CONCLUSIONS}

This work introduces a new class of generalized Mixed Poisson processes. First a new structure distribution is defined. It is called Min-U-Exp distribution. It is very similar to the exponential one and coincides with the distribution of the minima of two random variables - Uniform and Exponential. The inter-arrival times of these processes are described by newly-introduced Exp-Min-U-Exp distribution. The probability type of the moments of arrivals of the corresponding events is called Erlang-Min-U-Exp one, and its properties are thoroughly investigated. Along with our work some new multivariate distributions are defined. Ordered Mixed Poisson-Min-U-Exp distribution describes, for example, the joint distribution of the time-intersection of Mixed Poisson process with Min-U-Exp mixing variable. The corresponding distribution of the additive increments (which are dependent) is Mixed Poisson-Min-U-Exp one. The joint distribution of the inter-arrival times (which are also dependent) is Multivatiate Exp-Min-U-Exp distribution of $II^{-nd}$ kind. Different properties of these distributions are obtained, and their numerical characteristics are computed. The relations between the considered random elements are shown.

\bigskip
\section{ACKNOWLEDGMENTS} The work was supported by the Scientific Research Fund in Konstantin Preslavsky University of Shumen, Bulgaria under Grant Number RD-08-.../......2024 and project Number 2024 - FNSE – ...., financed by Scientific Research Fund of Ruse University.}

\bibliographystyle{plain}

\end{document}